\documentclass[final,3p,mathptmx]{elsarticle}

\usepackage{amsmath,amsthm,amssymb,amsfonts}
\usepackage{bm}
\usepackage{color}
\usepackage{rotating}
\usepackage{subcaption}
\usepackage{psfrag,graphicx,pstool}
\usepackage{varwidth}

\usepackage{tabularx}
\usepackage{ragged2e}
\newcolumntype{M}[1]{>{\centering\arraybackslash}m{#1}}
\usepackage{tikz}
\usetikzlibrary{plotmarks}
\usetikzlibrary{arrows}
\usetikzlibrary{shapes}

\newcommand{\lspan}{\mathop{\mathrm{span}}}

\theoremstyle{plain}
\newtheorem{thm}{Theorem}[section]
\newtheorem{lem}[thm]{Lemma}
\newtheorem{prop}[thm]{Proposition}

\theoremstyle{definition}
\newtheorem{defi}{Definition}[section]

\theoremstyle{remark}
\newtheorem{remark}[thm]{Remark}

\newcommand{\Wg}{E_{\Gamma}}

\usepackage{url}
\usepackage[
pdftex,
colorlinks=true,
draft=false,
bookmarks=false,
bookmarksnumbered=true,
plainpages=false,
]{hyperref}

\makeatletter
\renewcommand*\env@matrix[1][*\c@MaxMatrixCols c]{%
	\hskip -\arraycolsep
	\let\@ifnextchar\new@ifnextchar
	\array{#1}}
\makeatother

\biboptions{sort&compress}

\newcommand{\R}{\mathbb{R}}
\newcommand{\U}{\mathbb{U}}
\newcommand{\D}{\mathbb{D}}
\newcommand{\G}{\mathbb{G}}
\newcommand{\N}{\mathbb{N}}
\newcommand{\Z}{\mathbb{Z}}
\newcommand{\C}{\mathbb{C}}
\newcommand{\bv}{\mathbf{v}}
\newcommand{\bu}{\mathbf{u}}

\newcommand{\bw}{\mathbf{w}}
\newcommand{\bz}{\mathbf{z}}
\newcommand{\bff}{\mathbf{f}}
\newcommand{\bzero}{\mathbf{0}}

\newcommand{\lra}{\longrightarrow}

\newcommand{\cD}{\mathcal{D}}

\renewcommand{\a}{\alpha}

\newcommand\bgamma{{\boldsymbol\gamma}}
\newcommand\bmu{{\boldsymbol\mu}}
\newcommand\bal{{\boldsymbol\alpha}}

\begin{document}

\title{Annihilation Operators for Exponential Spaces in Subdivision}
\author[cc]{Costanza Conti\corref{cor1}}
\ead{costanza.conti@unifi.it}
\author[sl]{Sergio L\'opez-Ure\~na}
\ead{sergio.lopez-urena@uv.es}
\author[lr]{Lucia Romani}
\ead{lucia.romani@unibo.it}
\cortext[cor1]{Corresponding author}
\address[cc]{Dip. di Ingegneria Industriale, Universit\`a di Firenze, Viale Morgagni 40/44, 50134 Firenze, Italy}

\address[sl]{Dep. de Matem\`atiques, Universitat de Val\`encia, Doctor Moliner Street 50, 46100, Burjassot, Valencia, Spain}

\address[lr]{Dip. di Matematica, Alma Mater Studiorum Universit\`{a} di Bologna, Piazza P.S. Donato 5, 40126 Bologna, Italy}

\begin{abstract}
We investigate  properties of differential and difference operators annihilating certain finite-dimensional spaces of exponential functions in two variables that are connected to the representation of  real-valued trigonometric and hyperbolic functions.
Although exponential functions appear in a variety of contexts, the motivation behind this technical note comes from considering subdivision schemes where annihilation operators play an important role. Indeed, subdivision schemes with the capability of preserving exponential functions can be used to obtain an exact description of surfaces parametrized in terms of trigonometric and hyperbolic functions, and annihilation operators are useful to automatically detect the frequencies of such functions.
\end{abstract}


\maketitle

\section{Introduction}\label{sec:1}
In this technical note we investigate properties of differential and difference operators annihilating a particular type of spaces of exponential functions (for short exponential spaces). The motivation to study these operators comes from considering non-stationary subdivision schemes preserving exponential polynomials (see, e.g., \cite{CD20,CCS16, Noi,Yoon}), which  are efficient iterative algorithms for the definition of curves and surfaces.
Non-stationary subdivision schemes reproducing exponential polynomials, first proposed in \cite{NiraLuzzatto},
 have been at the center of a series of papers due to their importance in applications  like  geometric modelling (see, e.g., \cite{LVJIMA}) or image analysis (see, e.g., \cite{Unser}).
Up to now, in all subdivision schemes reproducing exponential polynomials, the exponential frequencies and their multiplicities are explicitly involved in the refinement rules. This is a major drawback since, usually, there is no a priori knowledge of them from the initial data.
Since a common machinery in designing linear, non-stationary subdivision schemes is the identification  of both a difference operator that annihilates the function space of interest, and a subdivision operator that reproduces it \cite{Tomas}, the aim of this technical note is twofold. On one hand, we investigate properties of differential and difference operators annihilating the type of exponential spaces of major interest in subdivision; on the other hand, following the seminal work of \cite{NiraLuzzatto}, we show how to use them for the automatic detection of the unknown exponential frequencies.
These ideas are already used in the univariate non-linear case (see \cite{SergioRosa}), where the derivation of subdivision rules that guarantee the preservation  of exponential polynomials rely on the definition of an \emph{annihilation operator} (also called \emph{annihilator}) whose kernel consists of them. In \cite{SergioRosa}, the construction of an interpolatory, non-linear, stationary subdivision scheme capable of reproducing functions in $\lspan\{1, \exp(\gamma z), \exp(-\gamma z)\}$,  where  $\gamma \in \G\setminus \{0\}$ and $\G:=\R_{\ge 0}\cup \imath (0,\pi)$, is indeed based on an annihilation operator.
For $f: \R\lra \C$,
this annihilator is obtained by the repeated application of the differential operator $D^{\gamma}f(z)=f'(z)-\gamma f(z),\ z\in \R.$ For the considered space, the annihilator is $D^{0}D^{\gamma}D^{-\gamma}f(z)=f'''(z)-\gamma^2 f'(z)$ since $$D^{0}D^{\gamma}D^{-\gamma}f(z)=0, \  \forall z \in \R \quad \Leftrightarrow \quad f\in \lspan\{1, \exp(\gamma z), \exp(-\gamma z) \, \} \quad \hbox{with}\quad \gamma \in \G\setminus \{0\}.$$
Replacing the differential operator $D^{\gamma}$ with the difference operator $\Delta^{\gamma}_t$ defined as   $\Delta^{\gamma}_t f(z)=f(z +t)-\exp(\gamma\, t)f(z),\ t,\, z\in \R$, the discrete operator $\Delta^0_t\Delta^{\gamma}_t \Delta^{-\gamma}_t$ turns out to be a discrete version of the annihilator $D^{0}D^{\gamma}D^{-\gamma}$ in the sense that
\begin{equation}\label{eq:diff}
\Delta^0_t\Delta^{\gamma}_t \Delta^{-\gamma}_t f(z)=0, \quad\forall t,\, z \in \R \quad \Leftrightarrow \quad f\in \lspan\{1, \exp(\gamma z), \exp(-\gamma z) \, \},\quad  \gamma \in \G\setminus \{0\}.
\end{equation}
If applied to a \emph{discrete} function $\bff^k:=\{f(2^{-k}\a),\  \a \in \Z\}$, $k \in \N$, basic ingredient of any subdivision scheme,
 the identity in \eqref{eq:diff}, for $t=2^{-k}$ and $z=2^{-k}(\a-1)$, yields an equation satisfied for all $\alpha\in \Z$ and
$f\in \lspan\{1, \exp(\gamma z), \exp(-\gamma z) \, \},\ \gamma \in  \G\setminus \{0\} $, that reads as
\begin{equation}\label{eq:rosasergio}
f(2^{-k}(\a-1))  -(2\cosh(2^{-k}\gamma)+1)f(2^{-k}\a)
+(2\cosh(2^{-k}\gamma)+1)f(2^{-k}(\a+1))-f(2^{-k}(\a+2))=0.
\end{equation}
Equation \eqref{eq:rosasergio} can be used to compute $\cosh(2^{-k}\gamma)$, that is to identify $\gamma$, and to set the
subdivision refinement rules able to reproduce the above-mentioned space of exponential functions (see \cite{SergioRosa} for all details).

The goal of this work is to investigate a similar idea for the bivariate case. We remark that we are aware that differential/difference operators associated with exponential functions are of interest in several domains of mathematics (other than subdivision schemes) and that we are omitting many related works. For example, \cite{Ron} provides a complete characterization of exponential polynomials using finite difference  operators in the very general but certainly more involved context of exponential Box-splines. But, here, our aim is to provide a simple and focused analysis of these operators helpful to understand the reproduction capabilities of bivariate subdivision schemes and related topics. Indeed, the results we here present on bivariate exponential spaces are complete and easy-to-follow, and make their understanding and use simple.

The rest of this paper  consists of two main sections. Section \ref{sec:characterization_differential} characterizes spaces of exponential functions as the kernel of an annihilation operator based on the repeated application of a particular differential operator. In Section \ref{sec:characterization_differential} we also show how this differential operator can be replaced by a discrete operator, which is the generalization to the bivariate setting of the finite difference  operator in \eqref{eq:diff}, proposed for univariate functions in \cite{SergioRosa}. Section \ref{sec:application} shows how to apply  annihilation operators  in the context of bivariate subdivision schemes that reproduce an important class of exponential functions. Conclusions are drawn in Section \ref{sec:conclusion}.

\section{Characterization of bivariate exponential functions via differential and difference operators} \label{sec:characterization_differential}

Let us start by defining the space of exponential functions we want to work with, relying on the standard column vector notation.
\begin{defi}\label{defi:exponential_polynomials}
Let $n \in \N$, $\D:=\R \cup \imath (-\pi,\pi)$,
\begin{equation}\label{def:gamma}
\bgamma^\ell \in \D^2, \, \ell=1, \ldots, n \quad \hbox{with} \quad \bgamma^\ell\neq \bgamma^j\quad \hbox{if}\quad \ell\neq j, \quad \hbox{and} \quad
\Gamma:=\left \{ \bgamma^1, \ldots, \bgamma^n \right \}.
\end{equation}
The space of exponential  functions associated to the set of $n$ distinct frequencies $\Gamma$ is
$$E_{\Gamma} := \lspan \left \{ \exp(\left(\bgamma^1\right)^T\bz ), \, \exp(\left(\bgamma^2\right)^T\bz ), \ldots, \exp(\left(\bgamma^n\right)^T\bz )  \ : \ \bgamma^\ell \in \Gamma, \, \ell=1,2,\ldots,n \right \},$$
where $\bz \in \R^2$ is the function variable.
\end{defi}

\begin{defi}\label{def:operator}
Let $F:\R^2\lra\C$ be a differentiable function and let $\nabla F(\bz)=\left(\frac{\partial F}{\partial z_1},\frac{\partial F}{\partial z_2}\right)^T$. Given $\bv\in\R^2 \setminus \{\bzero\}$ and $\bgamma\in\D^2$, we define the differential operator
\begin{equation}\label{eq:diffoper}
\cD_\bv^\bgamma F(\bz) := \left(\nabla F(\bz)-\bgamma F(\bz)\right)^T \bv, \qquad \bz\in \R^2.\end{equation}
\end{defi}

\noindent
Our aim is first, to show that the exponential functions in $E_{\Gamma}$ are in the kernel of the  differential operator $\cD_{\bv^1}^{\bgamma^1} \cD_{\bv^2}^{\bgamma^2} \ldots \cD_{\bv^n}^{\bgamma^n}$, where ${\bv^1},\ {\bv^2},\ \ldots, {\bv^n} \in \R^2 \setminus \{\bzero\}$ is an arbitrary set of $n$ non-zero directions. Second, to show that such exponential functions are also in the kernel of the discrete version of $\cD_{\bv^1}^{\bgamma^1} \cD_{\bv^2}^{\bgamma^2} \ldots \cD_{\bv^n}^{\bgamma^n}$.

\begin{remark}\label{rem0}
Let $\bv=(v_1,v_2)^T\in\R^2 \setminus \{\bzero\}$ and $|\bv|:=\sqrt{\bv^T \bv}=\sqrt{v_1^2+v_2^2}$.
It is easy to show that $\cD_\bv^\bgamma F(\bz) = 0 \Leftrightarrow  \cD_{\bv/|\bv|}^\bgamma F(\bz) = 0$, for all $\bz\in \R^2$.
Thus, from now on, we continue by assuming $\bv$ to be a unit vector of $\R^2$ and by writing $\bv\in \U:=\{ \bu \in \R^2 \, : \, |\bu|=1\}$.
\end{remark}

We go on with two preliminary results needed as basic steps of an induction argument used to prove Theorem \ref{prop:solution2}. The latter discusses solutions to differential equations of the type $\cD_{\bv^1}^{\bgamma^1} \cD_{\bv^2}^{\bgamma^2} \ldots \cD_{\bv^n}^{\bgamma^n} F=0$ where ${\bv^1},\ {\bv^2},\ \ldots, {\bv^n} \in \U$ is an arbitrary set of $n$ unit vectors, and  shows  that such solutions are independent of the vectors $\bv^1,\bv^2,\ldots,\bv^n$. \\

The first preliminary result considers the homogeneous case.

\begin{prop} \label{prop:solution_perp}
Let $\bv\in \U$ and $\bgamma \in\D^2$ be such that $\bgamma^T\bv \neq 0$ and $\bgamma^T\bv^\perp \neq 0$ where $\bv^\perp:=(-v_2,v_1)^T$. Then,
$$
\cD_\bw^\bgamma F= 0, \; \bw\in \{\bv, \bv^\perp\}\quad
\Leftrightarrow \quad F(\bz) = c \exp\left (\bgamma^T \bz \right ) \quad \hbox{with} \quad  c \in \C,\ \ \bz \in \R^2.
$$
\end{prop}
\proof
The implication $\Leftarrow$ is easily verified. To prove $\Rightarrow$ we know by classical arguments (see, e.g., \cite[chapter 3]{book}), applied for both $\bv$ and $\bv^\perp$, that
\begin{align*}
F(\bz) = \kappa_1\left (\bz^T \bv^\perp \right ) \exp\left ((\bgamma^T \bv) (\bz^T \bv) \right ) = \kappa_2\left (- \bz^T \bv \right ) \exp\left ((\bgamma^T \bv^\perp) (\bz^T \bv^\perp) \right ),
\end{align*}
where $\kappa_1,\ \kappa_2 :\R\lra\C$.
Taking $\bz=t\bv$, $t \in \R$, we obtain $\kappa_1 (0) \exp\left (\bgamma^T \bv t \right ) = \kappa_2\left (-t \right )$,
where we used that $\bv\in\U$. Denoting $c=\kappa_1(0)$, we can write (from the expression of $F(\bz)$ involving $\kappa_2$) that
$$
\begin{array}{c}
F(\bz)= c \exp\left ((\bgamma^T \bv) (\bz^T \bv) \right ) \exp\left ((\bgamma^T \bv^\perp) (\bz^T \bv^\perp) \right )
=c \exp\left ( (\bgamma^T \bv) (\bv^T \bz) \right ) \exp\left ((\bgamma^T \bv^\perp) ((\bv^\perp)^{T} \bz) \right ).
\end{array}
$$
Next we apply the additive property of exponentials and the distributivity of matrix operations, so obtaining
$$ F(\bz) =
c \exp \Big (\bgamma^T \left(\bv \bv^T +\bv^\perp (\bv^{\perp})^T \right)\bz  \Big ).
$$
Since $ \bv \bv^T +\bv^\perp (\bv^\perp)^{T} =
\begin{pmatrix}
1 & 0 \\ 0 & 1
\end{pmatrix}$,
we arrive at $F(\bz) = c \exp\left (\bgamma^T \bz  \right )$, so concluding the proof.
\qed

\smallskip
The second preliminary result considers the non-homogeneous case and can be proven by standard arguments (again look at \cite[chapter 3]{book}, for example).
\begin{lem} \label{prop:solution_complete_n}
Let  $\bv\in \U$, $\bgamma^1,\bgamma^2,\ldots,\bgamma^n$ be defined as in \eqref{def:gamma} and $(\bgamma^\ell-\bgamma^n)^T \bv \neq 0$, $\ell=1,\ldots,n-1$.
Then, for  $c_\ell \in \C$, $\ell =1, \ldots, n-1$ there exist $d_{\ell}(\bv) \in \C$, ${\ell}=1, \ldots, n-1$ such that
$$
\cD_\bv^{\bgamma^n} F(\bz) = \sum_{\ell=1}^{n-1} c_\ell \exp( \left( \bgamma^\ell \right)^T \bz)
\quad \Rightarrow \quad
F(\bz) =\sum_{{\ell}=1}^{n-1} d_{\ell}(\bv) \exp((\bgamma^{\ell})^T \bz)+ \kappa\left (\bz^T \bv^\perp \right ) \exp\left (((\bgamma^n)^T \bv) (\bz^T \bv) \right )
$$
where $\kappa:\R\lra\C$ is a particular function.
\end{lem}

Exploiting Lemma \ref{prop:solution_complete_n}, with an induction argument, we arrive at the sought result.

\begin{thm} \label{prop:solution2}
Given $\bgamma^1, \bgamma^2, \ldots, \bgamma^n$ as in \eqref{def:gamma} we have that,
$$
\begin{array}{c}
\cD_{\bv^1}^{\bgamma^1} \cD_{\bv^2}^{\bgamma^2} \ldots \cD_{\bv^n}^{\bgamma^n}  F= 0, \quad \forall \ \bv^1,\bv^2,\ldots,\bv^n \in\U \ \  \Leftrightarrow \ F(\bz) = \displaystyle \sum_{\ell=1}^n c_\ell \exp \left( (\bgamma^\ell)^T \bz \right) \ \hbox{with} \ c_\ell \in \C,\ \bz \in \R^2.
\end{array}
$$ 
\end{thm}

%
\proof
To show $\Leftarrow$ we first exploit the linearity and then the commutativity of the differential operator, so obtaining
$$
\cD_{\bv^1}^{\bgamma^1} \cD_{\bv^2}^{\bgamma^2} \ldots \cD_{\bv^n}^{\bgamma^n}  F(\bz)=
\sum_{\ell=1}^n c_\ell \, \cD_{\bv^1}^{\bgamma^1} \ldots \cD_{\bv^\ell}^{\bgamma^\ell} \ldots \cD_{\bv^n}^{\bgamma^n} \,  \exp((\bgamma^\ell)^T \bz)= \sum_{\ell=1}^n c_\ell \, \cD_{\bv^1}^{\bgamma^1} \ldots \cD_{\bv^n}^{\bgamma^n} \, \cD_{\bv^\ell}^{\bgamma^\ell} \exp((\bgamma^\ell)^T \bz) = 0.
$$
To prove $\Rightarrow$ we use induction on $n$. 
Proposition \ref{prop:solution_perp} proves the case $n=1$.  For a general $n$,
let $\bv$ be a unit vector such that $\bv^T (\bgamma^\ell-\bgamma^n)\neq 0$,
$\ell=1,2,\ldots n-1$.
Applying the induction hypothesis, we have
$$
 \cD_{\bv^1}^{\bgamma^1} \cD_{\bv^2}^{\bgamma^2} \ldots \cD_{\bv^{n-1}}^{\bgamma^{n-1}} \cD_{\bv}^{\bgamma^n}  F(\bz) = 0, \  \forall \bv^1,\bv^2,\ldots,\bv^{n-1} \in\U
\ \  \Rightarrow \ \
 \cD_{\bv}^{\bgamma^n} F(\bz)=   \sum_{\ell=1}^{n-1} c_{\ell}(\bv) \exp((\bgamma^\ell)^T \bz),
$$
where the dependence on $\bv$ of the coefficients  $c_{\ell}(\bv)$ is due to the fact that,  for each value of $\bv$, $\cD_{\bv}^{\bgamma^n} F$ is a different function.
 By Lemma \ref{prop:solution_complete_n} there exists  $\kappa:\R\lra\C$ such that
\begin{equation}\label{eq:sysF}
 F(\bz)=   \sum_{\ell=1}^{n-1} d_{\ell}(\bv) \exp((\bgamma^\ell)^T \bz) + \kappa(\bz^T \bv^\perp  )\exp\left (((\bgamma^n)^T \bv) (\bz^T \bv) \right ).
 \end{equation}
Setting $\bv^{1}=\bv^{2}=\ldots=\bv^{n-1}=\bv$, we have
\begin{equation}\label{eq:SERGIO}
\cD_{\bv}^{\bgamma^1} \cD_{\bv}^{\bgamma^2} \ldots \cD_{\bv}^{\bgamma^{n-1}} F(\bz)=
\kappa(\bz^T \bv^\perp  )\exp\left (((\bgamma^n)^T \bv) (\bz^T \bv) \right )\prod_{\ell=1}^{n-1} (\bgamma^\ell- \bgamma^n )^T \bv.
\end{equation}
Exploiting the commutativity property and using  the induction for $n=1$, we obtain
$$
 \cD_{\bw}^{\bgamma^n} \left (\cD_{\bv}^{\bgamma^1} \cD_{\bv}^{\bgamma^2} \ldots \cD_{\bv}^{\bgamma^{n-1}} F(\bz) \right )=   0, \  \forall \bw\in \U\quad
 \Longrightarrow \quad
\cD_{\bv}^{\bgamma^1} \cD_{\bv}^{\bgamma^2} \ldots \cD_{\bv}^{\bgamma^{n-1}} F(\bz) = C(\bv) \exp((\bgamma^n)^T \bz), 
$$
where $C(\bv) \in \C$ is again dependent on $\bv$ since  we apply the induction for each choice of a preliminary fixed $\bv$.
Comparing \eqref{eq:SERGIO} with the last expression of $\cD_{\bv}^{\bgamma^1} \cD_{\bv}^{\bgamma^2} \ldots \cD_{\bv}^{\bgamma^{n-1}} F(\bz)$ we find that,
$$
\kappa(\bz^T \bv^\perp  )\exp\left ( ((\bgamma^n)^T \bv) (\bz^T \bv) \right )
= c_n(\bv) \exp((\bgamma^n)^T \bz),\quad \hbox{for some}\quad  c_n(\bv)\in\C.
$$
Thus, replacing the latter in \eqref{eq:sysF}, we conclude that
$F(\bz)=   \sum_{\ell=1}^{n} c_\ell(\bv) \exp((\bgamma^\ell)^T \bz),$ with $c_\ell(\bv) =d_\ell(\bv),\ \ell=1, \ldots, n-1$.
To deduce that $c_\ell$, $\ell=1,\ldots,n$, do not depend on $\bv$,
we take an arbitrary $\bv$  inside the open set
$\{\bw \in \U \ : \ \bw^T (\bgamma^\ell-\bgamma^n)\neq 0,\ \ell=1,\ldots,n-1\}$,
and  differentiate $F$ with respect to $\bv$, providing
$0=   \sum_{\ell=1}^{n}  \left (\nabla c_\ell (\bv)\right )\exp((\bgamma^\ell)^T \bz).$
From the latter and the linear independence of $\{\exp((\bgamma^\ell)^T \bz)\}_{\ell=1}^n$, we
deduce that
$c_\ell,\ \ell=1, \ldots, n$, are indeed independent of $\bv$. \qed

\smallskip
The annihilation operator obtained via the differential operator in \eqref{eq:diffoper} admits a discrete analogue based on the repeated application of the difference operator defined next. Note that the presence of the scalar $t$ is to deal with non-unit vectors.

\begin{defi}
Given $F:\R^2\lra\C$, $\bv\in\U$, $\bgamma\in\D^2$,  we define the difference operator
$$\Delta_{t\bv}^\bgamma F(\bz) := F(\bz+t\bv)-\exp(\bgamma^T \bv t) F(\bz),\quad \bz\in \R^2,\quad t\in \R_{+}.$$
\end{defi}

\smallskip In the next theorem we show that $\cD^\bgamma_\bv$ can be replaced by $\Delta_{t \bv}^\bgamma$, $t \in \mathbb{R}_+$, to characterize functions in $E_{\Gamma}$.

\begin{thm} \label{thm:solution2}
Given $\bgamma^1,\bgamma^2,\ldots,\bgamma^n$ as in \eqref{def:gamma} we have that
$$\Delta_{t_1\bv^1}^{\bgamma^1} \Delta_{t_2\bv^{2}}^{\bgamma^{2}} \ldots \Delta_{t_n\bv^n}^{\bgamma^n} F= 0,
\ \ \forall\ t_i\bv^i\in\R^2\setminus \{{\bf 0}\},\ i=1,\ldots,n\ \ \  \Leftrightarrow \ F(\bz) = \displaystyle  \sum_{\ell=1}^n c_\ell \exp((\bgamma^\ell)^T \bz)  \ \hbox{with} \ c_\ell \in \C, \bz \in \R^2.$$
\end{thm}
\proof
We will prove that
	\begin{equation} \label{eq:delta_equations}
	\Delta_{t_1\bv^1}^{\bgamma^1} \Delta_{t_2\bv^{2}}^{\bgamma^{2}} \ldots \Delta_{t_n\bv^n}^{\bgamma^n} F= 0
\quad \hbox{if and only if}\quad
	\cD_{\bv^1}^{\bgamma^1} \cD_{\bv^2}^{\bgamma^2} \ldots \cD_{\bv^n}^{\bgamma^n}  F(\bz) = 0,\quad \ \ \forall\ t_i\bv^i \in\R^2\setminus \{{\bf 0}\},\ i=1,\ldots,n.
	\end{equation}
The left implication can be easily proven with the help of Theorem  \ref{prop:solution2}, while for the right implication we use an induction argument. For $n=1$,
we know that  $\Delta_{t\bv}^\bgamma F= 0$ means
$F(\bz+t\bv) = \exp(\bgamma^T \bv\ t)F(\bz)$ from which  we get
$$\nabla F(\bz)^T \bv = \lim_{t\rightarrow 0} t^{-1}(F(\bz+t\bv)-F(\bz)) =
F(\bz) \lim_{t\rightarrow 0} t^{-1}(\exp(\bgamma^T \bv t)-1) = F(\bz) \bgamma^T \bv,
$$ and therefore $\cD_\bv^\bgamma F(\bz) = (\nabla F(\bz) - \bgamma F(\bz))^T \bv = 0$.
For $n>1$, we use the easy-to-check commutativity of the differential and difference operators, i.e.,
$$\cD_\bv^\bgamma \Delta_{t\bw}^\mu = \Delta_{t\bw}^\mu \cD_\bv^\bgamma,\quad \forall \bgamma,\bmu\in \D^2,\ \bv,\bw\in\U,\quad t\in \R_{+}.$$
From above we see by induction that
for $\bv^1,\bv^2,\ldots,\bv^n\in\U$ and $\bgamma^1,\bgamma^2,\ldots,\bgamma^n$ as in \eqref{def:gamma},
\begin{equation}
 \cD_{\bv^1}^{\bgamma^1} \cD_{\bv^2}^{\bgamma^2} \ldots \cD_{\bv^n}^{\bgamma^n}  F= 0 \quad \Leftrightarrow \quad \Delta_{t_1\bv^1}^{\bgamma^1} \Delta_{t_2\bv^2}^{\bgamma^2} \ldots \Delta_{t_n\bv^n}^{\bgamma^n}  F= 0, \quad \forall t_i\in\R_+,\ i=1, \ldots, n,
 \end{equation}
 which concludes the proof.
  \qed

\begin{remark}
Note that,  when ${t\bv}\in \Z^2$, the operator  $\Delta_{t\bv}^\bgamma$ can also be applied to discrete data, as needed when we deal with subdivision schemes. Indeed, for the sequence ${\bf F}=\{F(\bal),\, \bal\in \Z^2\}$ seen as a function from $\Z^2$ to $\C$, the action of $\Delta_{t\bv}^\bgamma$ reads as
$$\Delta_{t\bv}^\bgamma F(\bal) := F(\bal+t\bv)-\exp(\bgamma^T \bv t) F(\bal),\quad  \bal \in\Z^2,\  t \bv \in\Z^2\setminus \{\bf 0\}.$$
Note also that, for a fixed $\bal$, the repeated application of the operator  $\Delta_{t\bv}^\bgamma$ with $t\bv\in \Z^2\setminus \{\bf 0\}$, namely $\Delta_{t_1\bv^1}^{\bgamma^1} \ldots \Delta_{t_n\bv^n}^{\bgamma^n} F(\bal)$, can be made local, i.e. involving only points of $\Z^2$ around the point $\bal$.
\end{remark}

\section{Applying annihilation operators in bivariate subdivision}\label{sec:application}
We continue this technical note with a section discussing  the use of bivariate, discrete annihilation operators that provide the natural generalization of the univariate, discrete annihilation operator in \eqref{eq:diff}. Such operators are applied to discrete data and used to identify a special space of exponential  functions of interest in subdivision. This is spanned by real-valued exponential functions obtained by assuming $n$ odd and considering the symmetric set of frequencies
$$
\Gamma=\{{\bf 0}\} \cup \{\bgamma^{\ell}, -\bgamma^{\ell}\}_{\ell=1,\ldots,\frac{n-1}{2}} \qquad \hbox{with} \quad \bgamma^{\ell} \in \G^2 \setminus \{{\bf 0}\}\quad \hbox{and} \quad \G:=\R_{\geq 0}\cup \imath (0,\pi).
$$
Specifically, in the following we focus our attention on $\Gamma=\{\bzero,\bgamma, -\bgamma, \tilde\bgamma, -\tilde \bgamma\}$ with $\bgamma:=(\gamma_1, \gamma_2)^T$, $\tilde \bgamma:=(\gamma_1, -\gamma_2)^T$ and $\bgamma, \tilde \bgamma \in \G^2 \setminus \{{\bf 0}\}$,
  since this set of frequencies is the one needed to reproduce spheres, hyperbolic paraboloids and other quadric surfaces. In fact, the coordinate components of the parametric representations of such surfaces are bivariate functions in
\begin{equation} \label{eq:gamma}
E_\Gamma = \lspan\{ 1, \exp( \bgamma^T \bz),  \exp(-\bgamma^T \bz), \exp(\tilde\bgamma^T \bz), \exp(-\tilde\bgamma^T \bz) \,\}, \quad\bgamma, \tilde \bgamma \in \G^2 \setminus \{{\bf 0}\} .
\end{equation}

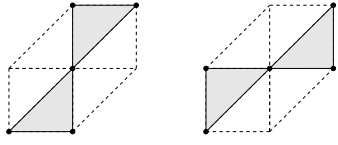
\begin{figure}[ht!]
\centering
\resizebox{!}{0.11\textwidth}{
\begin{tikzpicture}
		\draw[dashed] (-1,-1) -- (1,1);
		\draw[dashed] (-1,-1) -- (1,-1);
		\draw[dashed] (-1,-1) -- (-1,1);
		\draw[dashed] (-1,1) -- (1,3);
		\draw[dashed] (-1,1) -- (3,1);
		\draw[dashed] (1,-1) -- (3,1);
		\draw[dashed] (1,-1) -- (1,3);
		\draw[dashed] (1,3) -- (3,3);
		\draw[dashed] (3,1) -- (3,3);
		\draw[dashed] (1,1) -- (3,3);
		\draw[dashed] (-1,1) -- (1,1);
		\draw[dashed] (1,-1) -- (1,1);
        \filldraw[draw=black, fill=gray!20] (1,1) -- (3,3) -- (1,3) -- (1,1) -- cycle;
        \filldraw[draw=black, fill=gray!20] (1,1) -- (-1,-1) -- (1,-1) -- (1,1) -- cycle;
        \filldraw[black] (3,3) circle (2pt) node[anchor=south east] {};
        \filldraw[black] (1,3) circle (2pt) node[anchor=south east] {};
		\filldraw[black] (1,1) circle (2pt) node[anchor=south east] {};
		\filldraw[black] (1,-1) circle (2pt) node[anchor=south east] {};
		\filldraw[black] (-1,-1) circle (2pt) node[anchor=south east] {};
		\end{tikzpicture}}
\hspace{0.5cm}
\resizebox{!}{0.11\textwidth}{
\begin{tikzpicture}
		\draw[dashed] (-1,-1) -- (1,1);
		\draw[dashed] (-1,-1) -- (1,-1);
		\draw[dashed] (-1,-1) -- (-1,1);
		\draw[dashed] (-1,1) -- (1,3);
		\draw[dashed] (-1,1) -- (3,1);
		\draw[dashed] (1,-1) -- (3,1);
		\draw[dashed] (1,-1) -- (1,3);
		\draw[dashed] (1,3) -- (3,3);
		\draw[dashed] (3,1) -- (3,3);
		\draw[dashed] (1,1) -- (3,3);
		\draw[dashed] (-1,1) -- (1,1);
		\draw[dashed] (1,-1) -- (1,1);
        \filldraw[draw=black, fill=gray!20] (1,1) -- (-1,1) -- (-1,-1) -- (1,1) -- cycle;
        \filldraw[draw=black, fill=gray!20] (1,1) -- (3,1) -- (3,3) -- (1,1) -- cycle;
        \filldraw[black] (-1,-1) circle (2pt) node[anchor=south east] {};
        \filldraw[black] (-1,1) circle (2pt) node[anchor=south east] {};
		\filldraw[black] (1,1) circle (2pt) node[anchor=south east] {};
		\filldraw[black] (3,1) circle (2pt) node[anchor=south east] {};
		\filldraw[black] (3,3) circle (2pt) node[anchor=south east] {};
		\end{tikzpicture}
}
	\vspace{-0.05cm}
\caption{Pairs of triangles with vertices in $\Z^2$ used to check whether the function $F \in E_{\Gamma}$ is constant.}
	\label{fig:cosh_stencil}
\vspace{-0.4cm}
\end{figure}

\begin{remark}\label{rem:constant}
For later use, we observe that if $E_\Gamma$ is the space of  exponential functions  in \eqref{eq:gamma}, and a bivariate function $F \in E_\Gamma$ assumes a constant value at 5 distinct points of $\Z^2$ that correspond to the vertices of the pair of grey triangles displayed in each illustration of Figure  \ref{fig:cosh_stencil}, then $F$ is a constant function that can be associated to the frequency $\bgamma=\bf 0$.
\end{remark}

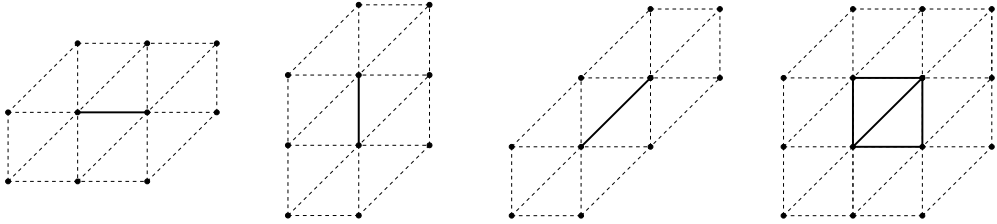
\begin{figure}[ht!]
\hspace{1.5cm}
	\resizebox{\textwidth}{!}{\begin{tabularx}{\textwidth}{M{0.18\textwidth} M{0.18\textwidth} M{0.18\textwidth} M{0.18\textwidth}}
			\resizebox{0.18\textwidth}{!}{
\begin{tikzpicture}
				\draw[ultra thick] (-1,0) -- (1,0);
				\draw[dashed] (-1,2) -- (3,2);
				\draw[dashed] (-3,-2) -- (1,-2);
				\draw[dashed] (-1,2) -- (-1,-2);
				\draw[dashed] (1,-2) -- (1,2);
				\draw[dashed] (-3,-2) -- (1,2);
				\draw[dashed] (-1,-2) -- (3,2);
				\draw[dashed] (-3,-2) -- (-3,0);
				\draw[dashed] (-1,0) -- (-3,0);
				\draw[dashed] (-1,2) -- (-3,0);
				\draw[dashed] (1,-2) -- (3,0);
				\draw[dashed] (1,0) -- (3,0);
				\draw[dashed] (3,2) -- (3,0);
				\filldraw[black] (-1,0) circle (2pt) node[anchor=east] {};
				\filldraw[black] (1,0) circle (2pt) node[anchor=west] {};
				\filldraw[black] (-1,2) circle (2pt) node[anchor=south] {};
				\filldraw[black] (1,2) circle (2pt) node[anchor=south] {};
				\filldraw[black] (3,2) circle (2pt) node[anchor=south] {};
				\filldraw[black] (-3,-2) circle (2pt) node[anchor=north] {};
				\filldraw[black] (-1,-2) circle (2pt) node[anchor=north] {};
				\filldraw[black] (1,-2) circle (2pt) node[anchor=north] {};
				\filldraw[black] (-3,0) circle (2pt) node[anchor=north] {};
				\filldraw[black] (3,0) circle (2pt) node[anchor=north] {};
			\end{tikzpicture}
} & \hspace{-0.4cm}
			\resizebox{!}{0.18\textwidth}{
\begin{tikzpicture}
				\draw[dashed] (-2,-3) -- (-2,1);
				\draw[ultra thick] (0,-1) -- (0,1);
				\draw[dashed] (2,-1) -- (2,3);
				\draw[dashed] (-2,-1) -- (2,-1);
				\draw[dashed] (-2,1) -- (2,1);
				\draw[dashed] (-2,-1) -- (2,3);
				\draw[dashed] (-2,-3) -- (2,1);
				\draw[dashed] (-2,-3) -- (0,-3);
				\draw[dashed] (0,-1) -- (0,-3);
				\draw[dashed] (2,-1) -- (0,-3);
				\draw[dashed] (-2,1) -- (0,3);
				\draw[dashed] (0,1) -- (0,3);
				\draw[dashed] (2,3) -- (0,3);
				\filldraw[black] (0,1) circle (2pt) node[anchor=south] {};
				\filldraw[black] (0,-1) circle (2pt) node[anchor=north] {};
				\filldraw[black] (2,3) circle (2pt) node[anchor=west] {};
				\filldraw[black] (2,1) circle (2pt) node[anchor=west] {};
				\filldraw[black] (2,-1) circle (2pt) node[anchor=west] {};
				\filldraw[black] (-2,1) circle (2pt) node[anchor=east] {};
				\filldraw[black] (-2,-1) circle (2pt) node[anchor=east] {};
				\filldraw[black] (-2,-3) circle (2pt) node[anchor=east] {};
				\filldraw[black] (0,-3) circle (2pt) node[anchor=north] {};
				\filldraw[black] (0,3) circle (2pt) node[anchor=north] {};
			\end{tikzpicture}
} & \hspace{-0.4cm}
			\resizebox{!}{0.18\textwidth}{
\begin{tikzpicture}
				\draw[ultra thick] (-1,-1) -- (1,1);
				\draw[dashed] (-3,-1) -- (1,3);
				\draw[dashed] (-1,-3) -- (3,1);
				\draw[dashed] (-3,-1) -- (1,-1);
				\draw[dashed] (-1,1) -- (3,1);
				\draw[dashed] (-1,-3) -- (-1,1);
				\draw[dashed] (1,-1) -- (1,3);
				\draw[dashed] (-1,-3) -- (-3,-3);
				\draw[dashed] (-3,-1) -- (-3,-3);
				\draw[dashed] (-1,-1) -- (-3,-3);
				\draw[dashed] (1,3) -- (3,3);
				\draw[dashed] (3,1) -- (3,3);
				\draw[dashed] (1,1) -- (3,3);
				\filldraw[black] (1,1) circle (2pt) node[anchor=north east] {};
				\filldraw[black] (1,1) circle (2pt) node[anchor=south west] {};
				\filldraw[black] (-3,-1) circle (2pt) node[anchor=south east] {};
				\filldraw[black] (-1,-3) circle (2pt) node[anchor=north west] {};
				\filldraw[black] (1,3) circle (2pt) node[anchor=south east] {};
				\filldraw[black] (3,1) circle (2pt) node[anchor=north west] {};
				\filldraw[black] (1,-1) circle (2pt) node[anchor=north west] {};
				\filldraw[black] (-1,1) circle (2pt) node[anchor=south east] {};
				\filldraw[black] (-1,-1) circle (2pt) node[anchor=south east] {};
				\filldraw[black] (3,3) circle (2pt) node[anchor=north west] {};
				\filldraw[black] (-3,-3) circle (2pt) node[anchor=south east] {};
			\end{tikzpicture}
} &
		\resizebox{!}{0.18\textwidth}{
\begin{tikzpicture}
		\draw[ultra thick] (-1,-1) -- (1,1);
		\draw[ultra thick] (-1,-1) -- (1,-1);
		\draw[ultra thick] (-1,-1) -- (-1,1);
		\draw[dashed] (-3,-1) -- (1,3);
		\draw[dashed] (-1,-3) -- (3,1);
		\draw[dashed] (-3,-1) -- (1,-1);
		\draw[dashed] (-1,1) -- (3,1);
		\draw[dashed] (-1,-3) -- (-1,1);
		\draw[dashed] (1,-1) -- (1,3);
		\draw[dashed] (-1,-3) -- (-3,-3);
		\draw[dashed] (-3,-1) -- (-3,-3);
		\draw[dashed] (-1,-1) -- (-3,-3);
		\draw[dashed] (1,3) -- (3,3);
		\draw[dashed] (3,1) -- (3,3);
		\draw[dashed] (1,1) -- (3,3);
		\draw[ultra thick] (-1,1) -- (1,1);
		\draw[dashed] (-1,3) -- (-1,-3);
		\draw[ultra thick] (1,-1) -- (1,1);
		\draw[dashed] (-3,-1) -- (-3,1);
		\draw[dashed] (-1,1) -- (-3,1);
		\draw[dashed] (-1,3) -- (-3,1);
		\draw[dashed] (3,3) -- (3,1);
		\draw[dashed] (1,3) -- (-1,3);
		\draw[dashed] (-1,-3) -- (1,-3);
		\draw[dashed] (1,-1) -- (1,-3);
		\draw[dashed] (1,-1) -- (3,-1);
		\draw[dashed] (3,1) -- (3,-1);
		\draw[dashed] (1,-3) -- (3,-1);
		\filldraw[black] (1,1) circle (2pt) node[anchor=north east] {};
		\filldraw[black] (1,1) circle (2pt) node[anchor=south west] {};
		\filldraw[black] (-3,-1) circle (2pt) node[anchor=south east] {};
		\filldraw[black] (-1,-3) circle (2pt) node[anchor=north west] {};
		\filldraw[black] (1,3) circle (2pt) node[anchor=south east] {};
		\filldraw[black] (3,1) circle (2pt) node[anchor=north west] {};
		\filldraw[black] (1,-1) circle (2pt) node[anchor=north west] {};
		\filldraw[black] (-1,1) circle (2pt) node[anchor=south east] {};
		\filldraw[black] (-1,-1) circle (2pt) node[anchor=south east] {};
		\filldraw[black] (3,3) circle (2pt) node[anchor=north west] {};
		\filldraw[black] (-3,-3) circle (2pt) node[anchor=south east] {};
		\filldraw[black] (1,-3) circle (2pt) node[anchor=south east] {};
		\filldraw[black] (3,-1) circle (2pt) node[anchor=south east] {};
		\filldraw[black] (-3,1) circle (2pt) node[anchor=south east] {};
		\filldraw[black] (-1,3) circle (2pt) node[anchor=south east] {};
		\end{tikzpicture}
}
	\end{tabularx}}
	\vspace{-0.1cm}
	\caption{From left to right: the horizontal, vertical and diagonal stencils of the extended butterfly subdivision scheme and the union of their points. The ticker edges correspond to
	 locations where the new points will be inserted by the refinement rules.}
	\label{fig:B_stencil}
\end{figure}

We have already mentioned in the Introduction that, to model parametric surfaces like spheres, hyperbolic paraboloids and quadric surfaces, we can
conveniently use subdivision schemes. Roughly speaking, subdivision schemes are efficient iterative methods capable of generating surfaces from samples of  a function on $\Z^2$ by means of the repeated application of simple and local refinement rules.  At level $k$, such rules use values on the grid $2^{-k}\Z^2$ to generate values on the grid  $2^{-k-1}\Z^2$, and so on (see for example the recent survey \cite{CD20}).
As it is well known, the refinement rules of subdivision schemes reproducing exponential polynomials are not only level-dependent, but  also explicitly involving the exponential frequencies $\bgamma$ (see, e.g., \cite{Noi,SergioRosa,NiraLuzzatto,Yoon,LVJIMA}). This is a major drawback since, usually, there is no a priori knowledge of $\bgamma$ from the initial data. Hence, the development of techniques to automatically estimate $\bgamma$ from the given data are very important.  The latter observation motivates a discussion on how discrete difference operators can be used as \emph{annihilation operators} for the automatic detection of exponential frequencies. 
In the context of subdivision schemes it is also important to use annihilators that are \emph{local} and  
involve the same sets of points -generally called \emph{stencils}-  used by the refinement rules of the scheme.
To keep the exposition simple, we stay focused both on the space \eqref{eq:gamma} (where $\Gamma=\{\bzero, \bgamma, -\bgamma, \tilde\bgamma, -\tilde \bgamma\}$ and whose importance is clear from the above discussion), and on the stencils of the \emph{extended butterfly} subdivision scheme (see \cite{RNY16} and references therein).  The  refinement rules of the extended  butterfly subdivision scheme  involve points
distributed as in Figure \ref{fig:B_stencil}.

From Theorem \ref{thm:solution2} we deduce that, for $F\in\Wg$ and $\mathbf{\bal} \in \mathbb{Z}^2$,
\begin{equation} \label{eq:Deltas}
\Delta^{0}_{{t_0\bv}^0}\; \Delta^{\boldsymbol{\gamma}}_{t_1\bv^1}\;\Delta^{-\boldsymbol{\gamma}}_{t_2\mathbf{v}^2}\;\Delta^{\widetilde{\boldsymbol{\gamma}}}_{t_3\mathbf{v}^3}\;\Delta^{-\widetilde{\boldsymbol{\gamma}}}_{t_4\mathbf{v}^4}\;F(\bal)\;=\;0,\quad \forall \bv^i\in\U, \ t_i\in \R_+ \ \hbox{such that} \ t_i\bv^i\in \Z^2,\ i=0,\ldots,4.
\end{equation}
Therefore, the operator $\mathcal{N}_{\bal,t_0\mathbf{v}^0,t_1\mathbf{v}^1,t_2\mathbf{v}^2,t_3\mathbf{v}^3,t_4\mathbf{v}^4}(F,\bgamma) = \Delta^{0}_{{t_0\bv}^0}\; \Delta^{\boldsymbol{\gamma}}_{t_1\bv^1}\;\Delta^{-\boldsymbol{\gamma}}_{t_2\mathbf{v}^2}\;\Delta^{\widetilde{\boldsymbol{\gamma}}}_{t_3\mathbf{v}^3}\;\Delta^{-\widetilde{\boldsymbol{\gamma}}}_{t_4\mathbf{v}^4}\;F(\bal)$
is an annihilation operator for $E_{\Gamma}$. Now, observe that for specific choices of $t_0\mathbf{v}^0, t_1\mathbf{v}^1, t_2\mathbf{v}^2, t_3\mathbf{v}^3, t_4\mathbf{v}^4$, some finite differences coincide.
For instance,
\begin{equation}\label{eq:symmetries}
\Delta^{\boldsymbol{\gamma}}_{(1,0)^T}\;=\;\Delta^{\widetilde{\boldsymbol{\gamma}}}_{(1,0)^T},\qquad\Delta^{\boldsymbol{\gamma}}_{(0,1)^T}\;=\;\Delta^{-\widetilde{\boldsymbol{\gamma}}}_{(0,1)^T},\qquad\Delta^{-\boldsymbol{\gamma}}_{(1,0)^T}\;=\;\Delta^{-\widetilde{\boldsymbol{\gamma}}}_{(1,0)^T}\quad\text{ and }\quad\Delta^{-\boldsymbol{\gamma}}_{(0,1)^T}\;=\;\Delta^{\widetilde{\boldsymbol{\gamma}}}_{(0,1)^T}.
\end{equation}
Hence, it turns out that to annihilate $E_\Gamma$ it is sufficient to consider just 3 finite difference operators, that is the annihilators
\begin{align} \label{eq:3_deltas}
\mathcal{N}_{\bal,t\mathbf{v},\mathbf{e},\mathbf{e}}(F,\bgamma)\; &=
\Delta^{0}_{t\mathbf{v}}\;\Delta^{\boldsymbol{\gamma}}_{\mathbf{e}}\;\Delta^{-\boldsymbol{\gamma}}_{\mathbf{e}}\;F(\bal), \quad \mathbf{e}\in\{(1,0)^T,(0,1)^T\}, \quad  \bal\in \Z^2,\  t\mathbf{v} \in \Z^2\setminus \{\bf 0\}.
\end{align}
Note that $\mathcal{N}_{\bal,t\mathbf{v},\mathbf{e},\mathbf{e}}(F,\bgamma)$ involves  function evaluations at
$\bal + \lambda \, t\mathbf{v}+ \mu \mathbf{e}$ with $\lambda \in \{0,1\}, \ \mu \in  \{0,1,2\}$, i.e. at 6 points of $\Z^2$ at most, and this is certainly a great advantage  compared with the 14 function evaluations that $\mathcal{N}_{\bal,t_0\mathbf{v}^0,t_1\mathbf{v}^1,t_2\mathbf{v}^2,t_3\mathbf{v}^3,t_4\mathbf{v}^4}(F,\bgamma)$ requires.
For the application of annihilators in the context of subdivision schemes this is a crucial feature.
Now, expansion of $\mathcal{N}_{\bal,t\mathbf{v},\mathbf{e},\mathbf{e}}(F,\bgamma)=0$ provides the equation
\begin{equation} \label{eq:explicit_annihilator}
-2\cosh(\boldsymbol{\gamma}^T\mathbf{e})\;\Delta^{0}_{t\mathbf{v}}F(\bal+\mathbf{e})\;+\;\Delta^{0}_{t\mathbf{v}}\;F(\bal+2\mathbf{e})\;+\;\Delta^{0}_{t\mathbf{v}}\;F(\bal)=0,\quad \alpha\in \Z^2,
\end{equation}
from which, whenever  $F\in E_\Gamma$ and $\Delta^{0}_{t\mathbf{v}}F(\bal+\mathbf{e})\neq0$, we can identify the unknown parameter $\boldsymbol{\gamma}$.
As a matter of fact, from \eqref{eq:explicit_annihilator} we can write
\begin{equation} \label{eq:cosh}
\cosh(\boldsymbol{\gamma}^T\mathbf{e})\;=\;\frac{\Delta^{0}_{t\mathbf{v}}F(\bal+2\mathbf{e})\;+\;\Delta^{0}_{t\mathbf{v}}F(\bal)}{2\;\Delta^{0}_{t\mathbf{v}}F(\bal+\mathbf{e})}, \qquad \forall \bal\in \mathbb{Z}^2,
\end{equation}
and, with the choices $\mathbf{e} = (1,0)^T$ and $\mathbf{e} = (0,1)^T$, we indeed univocally identify $\cosh(\gamma_1)$ and  $\cosh(\gamma_2)$, respectively, and hence $\bgamma$.
Note that, in case $\Delta^{0}_{t\mathbf{v}}F(\bal+\mathbf{e})=0$,  to make it non-zero we change  $t\mathbf{v}$ with the restriction of involving only points belonging to the union of the three subdivision stencils of the extended butterfly scheme shown in the rightmost picture of Figure \ref{fig:B_stencil}. According to Remark \ref{rem:constant}, for $\mathbf{e}=(1,0)^T$ and a fixed selection of $\bal$, there are two different sets of choices for $t\mathbf{v}$ that are enough to be considered (see Figure \ref{fig:cosh_stencil} where the vertex shared by the two triangles represents $\bal$). Precisely, as to the leftmost picture, the vectors are $t\mathbf{v} \in \{ (0,1)^T, \, (1,1)^T, \, (0,-1)^T, \, (-1,-1)^T \}$ whereas, as to the rightmost picture, they are
$t\mathbf{v} \in \{ (1,0)^T, \, (1,1)^T, \, (-1,0)^T, \, (-1,-1)^T \}$, respectively.
Therefore, if for a fixed selection of $\bal$, $\Delta^{0}_{t\mathbf{v}}F(\bal+\mathbf{e})=0$ for all the choices of $t\mathbf{v}$ in one of the  two sets mentioned above, this means that $F$ is constant, which anyway implies the identification of $\bgamma$ with ${\bf 0}$. The same obviously  holds for $\mathbf{e}=(0,1)^T$.

	\vspace{-0.1cm}
\section{Closing remarks}\label{sec:conclusion}
This technical note investigates properties of differential and difference operators annihilating certain finite-dimensional spaces of  bivariate exponential functions connected to the representation of  real-valued trigonometric and hyperbolic functions used in the context of subdivision schemes.
The example illustrated in Section \ref{sec:application} has shown the benefits of such annihilation operators to identify the frequency $\bgamma$ involved in the definition of subdivision rules. For further details concerning the use of  annihilators in subdivision schemes, we refer the interested reader to the follow-up of this theoretical work recently appeared in \cite{SA}. In that paper, annihilation operators as those in \eqref{eq:3_deltas}
are shown in action,  not only to detect frequencies, but also to locally identify the correct subdivision rule to be applied.
This demonstrates the important role that annihilators can play.
 
\bigskip {\bf Acknowledgement:} The first and the last author are members of INdAM-GNCS, Italy, partially supporting this work.

\end{document}